\documentclass[11pt,a4paper]{article}
\usepackage[centertags]{amsmath}
\usepackage{latexsym}
\usepackage{amsfonts}
\usepackage{amssymb}
\usepackage{amsthm}

\usepackage[english]{babel}

\newtheorem{remark}{Remark}

\theoremstyle{plain}
\newtheorem{lemma}{Lemma}
\newtheorem{proposition}{Proposition}
\newtheorem{theorem}{Theorem}
\newtheorem{definition}{Definition}

\newtheorem{example}{Example}

\title{\LARGE Structure of Symplectic Lie groups and momentum map}
\author{Alberto Medina \footnote{\footnotesize {Universit\'e de Montpellier 2, D\'epartement de Math\'ematiques, 34095, France. e-mail~: medina@math.univ-montp2.fr.}}}

\begin{document}
\maketitle

\begin{abstract}
We describe the structure of the Lie groups endowed with a left-invariant
symplectic form, called symplectic Lie groups, in terms of semi-direct products
of Lie groups, symplectic reduction and principal bundles with affine fiber.This description is particularly nice if the group is Hamiltonian, that is, if the left canonical action of the
group on itself is Hamiltonian. The principal tool used for our description is a canonical
affine structure associated with the symplectic form. We also characterize the Hamiltonian symplectic Lie groups among the connected symplectic Lie groups. We specialize our principal results to the cases of simply connected Hamiltonian symplectic nilpotent Lie groups or Frobenius symplectic Lie groups. Finally we pursue the study of the classical affine Lie group as a symplectic Lie group.
\end{abstract}

MSC Classes 53D20,70G65


Symplectic reduction,Symplectic double extension.
\section{Introduction}
A closed nondegenerate smooth 2-form ${\Omega}$ on a manifold $M$ is called a symplectic form and the pair
$(M,\Omega )$ is called a symplectic manifold. A smooth action $\Psi$ of a Lie group $G$ on $(M,\Omega )$
is symplectic if ${\Omega}$ is $\Psi$-invariant. This means that $i(\overline{x}){\Omega}$ is a closed 1-form on
 $M$ for any $x$ in the Lie algebra  ${\cal G}:=T_{\epsilon}(G)=L(G)$  of $G$, where $\overline{x}$ denotes the fundamental vector
 field associated with $x$ and $\epsilon$ is the unit of $G$. When $i(\overline{x}){\Omega}$ is exact on $M$ for any $x\in{\cal G}$, the action
 $\Psi$ is Hamiltonian i.e. there exists a smooth map $J:M\longrightarrow {\cal G}^{*},p\mapsto J(p)$ such that $i(\overline{x}){\Omega}$ is the derivative of the smooth map  $J_{x}: M\longrightarrow \mathbb R$ associated with each $x\in \cal G$ by $p\mapsto J_{x}(p):=J(p)(x)$.
 Such a map is called a momentum map for $\Psi$.

\begin{definition}
 A Lie group $G$ endowed with a left-invariant symplectic form $\omega ^{+}$
 is called a symplectic Lie group. The pair $(G,\omega^{+} )$ is  Hamiltonian
if the left canonical action $L_{G}$ of $G$ on $G$ (i.e.by means of left translations) is Hamiltonian.
The pair $({\cal G},\omega )$ with $\omega =\omega^{+}_\epsilon $, where $\epsilon$ is
the identity of $G$, is called a symplectic Lie algebra.
\end{definition}

\begin{definition} A pair $(G,\nabla^{+} )$, where ${\nabla^{+}}$ is a left-invariant affine
 structure on $G$, will be called an affine Lie group.
 \end{definition}
To a left-invariant affine structure on $G$ corresponds  a bilinear product $(x,y)\mapsto xy$ on ${\cal G}$ such that
\begin{equation}\label{eq:lsa}
    (xy)z-x(yz)=(yx)z-y(xz)
\end{equation}
This product is compatible with the Lie bracket of ${\cal G}$ in the sense that
\begin{center}
$xy-yx=[x,y]$
\end{center}
for any $x,y,z\in {\cal G}$.
 
A vector space $A$ endowed with a bilinear product that satisfies (1) is called 
a left-symmetric algebra. In this case $[x,y]:=xy-yx$ is a Lie bracket on A and the derived ideal $[A,A]$ is different from $A$ (\cite{He}).
\begin{remark}\label{rem1}
1.If $G$ admits a left-invariant symplectic form its Lie algebra is different
from its derived ideal $[\cal G,\cal G]$. In particular a semisimple Lie group is not a
symplectic Lie group.
\\
2. A symplectic Lie group $(G,\omega^{+} )$ is endowed with
a left-invariant affine structure given by the left-invariant torsion-free and flat connection
$\nabla^{+} $ defined by the formula:
\begin{equation}
    \omega^{+}(\nabla _{x}{+}y^{+},z^{+})=-\omega^{+}(y^{+},[x^{+},z^{+}])
\label{eq:ED1}
\end{equation}
where $a^{+}$ is the left-invariant vector field on $G$ associated with $a\in {\cal G}=T{\epsilon} (G)$.

\end{remark}
The aim of this paper is to describe the structure of  symplectic Lie groups in
terms of semi-direct products of Lie groups or symplectic reduction and principal
fiber bundles with affine fibers improving  results given in \cite{Da-Me2}.This description
can be bettered if the canonical left action $L_{G}$ of $G$ on $G$ is Hamiltonian.

The affine structure associated with $\omega^{+}$ plays a very important role in our
study.
Obviously a better understanding of the structure of symplectic Lie groups
is important in Mathematics and Physics.

\section{Structure of Symplectic Lie groups}
In this section we will describe the general structure of the symplectic Lie groups by means
of  Theorems \ref{theo2.3} and \ref{theo2.4} which complete and improve results given in \cite{Da-Me2}.

Let $(G,\omega^{+})$ be a connected symplectic Lie group. Consider a Lie subgroup $H$ of $G$ with Lie algebra $L(H)$and denote by $L_{H}$, respectively $R_{H}$, the left (resp. right) canonical action of $H$ on $G$ given by the product of $G$. Obviously $L_{H}$ is symplectic.

\medskip 

\noindent
{\bf Study of $L_{H}$ and $R_{H}$}.

Let $O_{L}(\sigma)=H{\sigma}$ be the orbit of ${\sigma}\in G $ and let $F$ be the subbundle of
$TG$, tangent to the orbits of $L_{H}$. Denote by ${\cal F}$
the foliation associated with $F$. It is well known that the subbundle
$F^{\bot }$, orthogonal to $F$, is integrable. Let ${\cal F^{{\bot }}}$ be
the corresponding foliation. The fiber of $F$ above ${\sigma}\in G$ is given by
$F_{\sigma}= \{x^{-}_{\sigma}; x\in L(H)\}$ with  $x^{-}$ the right-invariant vector field on $G$ associated with $x\in L(H)$. The foliation ${\cal F^{{\bot}}}$
is defined by the intersection of the kernels of the closed
1-forms $i(x^{-})\omega^{+}$ where $x\in L(H)$.
On the other hand the orbits of $R_{H}$ are the leaves of the integrable distribution
$D_{\sigma}^{+}:=\{x_{\sigma}^{+};x\in L(H)\}$.
For ${\tau}\in H$ we have:
\begin{center}
     $O_{L}(\tau)=H=O_{R}(\tau)$
\end{center}
From the formula $\omega_{\sigma}^{+}(x_{\sigma}^{+},y_{\sigma}^{+})=\omega(x,y)$ for $x,y \in {\cal G}$,
it follows that if $L(H)$ is symplectic (resp. totally isotropic or coisotropic or Lagrangian) then $R_{H}$ determines a left-invariant symplectic (resp. totally isotropic or coisotropic or Lagrangian) foliation on $G$.
Moreover if $H$ is normal in $G$ we obtain:

\begin{lemma}\label{lem2.1} Let $H$ be a connected normal Lie subgroup of a connected symplectic Lie group $(G,\omega^{+})$. Then ${\cal F^{{\bot}}}$ is left invariant and its leaves are affine submanifolds of $(G,\nabla)$. Moreover if $L_{H}$ is Hamiltonian, the leaves of ${\cal F^{{\bot}}}$ are closed.\end{lemma}

\proof For $a,b\in L(H)^{\bot}$ and $x\in L(H)$, we have
\begin{center}
$\omega(ab,x)=-\omega(b,[a,x])=0$.
\end{center}
Consequently $L(H)^{\bot}$ is  a subalgebra of the left-symmetric algebra $L(G)$ given by (\ref{eq:ED1}).

Suppose that $L_{H}$ is hamiltonian and let $J:G\to L(H)^{*}$ be a momentum map for $L_{H}$.
As the orbits of $L_{H}$ have the same dimension, the rank of $J$ is constant on $G$ and, for every $\sigma \in G$, the connected component of $J^{-1}(J(\sigma))$ containing $\sigma$ is the leaf of ${\cal F^{{\bot}}}$ through $\sigma$. Hence the leaves of ${\cal F^{{\bot}}}$ are closed.\qed

\begin{lemma}\label{lem2.2} Every symplectic Lie group $(G,\omega^{+} )$ contains a non-trivial connected normal and closed subgroup $H$ isomorphic to a vector group or a torus.\end{lemma}

\proof Since $G$ is an affine group, the Lie algebra ${\cal G}$ is different from its derived ideal $[{\cal G},{\cal G}]$ (see \cite{He}). Consequently
$G$ contains a subgroup verifying the conditions required in the Lemma (see for example \cite{He}).\qed

\begin{theorem}\label{theo2.3} Let $(G,\omega^{+} )$ be a connected symplectic Lie group and $H$ a non-trivial
connected normal and closed subgroup of $G$. Denote by $H^{{\bot}}$ the connected Lie subgroup
of $G$ orthogonal to $H$ and put $K=H\cap H^{{\bot}}$. We have

1. If $K$ is discrete, $H$ and $H^{{\bot}}$ are symplectic Lie subgroups of $(G,\omega^{+} )$
and $G$ contains an open connected  symplectic Lie subgroup isomorphic to group
$(H\times H^{{\bot}})/Ker f$ where $H\times H^{{\bot}}$ is a semi-direct product, $H^{{\bot}}$
acts on $H$ by symplectomorphisms and $f$ is the homomorphism of Lie groups from $H\times H^{{\bot}}$ into $G$ given by $f(\sigma ,\tau)=\sigma\tau$.

2. If $K$ is not discrete then it is a closed Abelian and normal subgroup of $H^{{\bot}}$,
the quotient group $R= H^{{\bot}}/K $ is a reduced symplectic Lie group of $(G,\omega^{+} )$
and the canonical sequence
\begin{equation}
\{\epsilon\} \rightarrow K\rightarrow H^{{\bot}}\rightarrow R\rightarrow \{\epsilon\}
\label{eq:ED3}
\end{equation}
is an exact sequence of affine Lie groups where the affine structures are deduced from
those of $(G,\nabla^{+} ),$ that of  $K$ being trivial. \end{theorem}

\proof The Lie algebra I of $H$ is a Lie ideal of ${\cal G}$. The identity
\begin{center}
$\omega (ab,x)=-\omega (b,[a,x])=0$
\end{center}
for $a,b\in I^{{\bot}}$ and $x\in I $, shows that $I^{{\bot}}$ is a left-symmetric
subalgebra of ${\cal G}$. Hence $H^{{\bot}}$ is an affine Lie subgroup of $(G,\nabla^{+})$.

Suppose $K$ is discrete. In this case ${\cal G}=I\times I^{{\bot}}$(semi-direct product). Obviously, the connected subgroup $K$ is closed in $H^{{\bot}}$. Moreover, for $a\in I\cap I^{{\bot}}, b\in I^{{\bot}}$ and $i\in I$, the identity
\begin{center}
$\omega ([a,b],i)+\omega ([b,i],a)+\omega ([i,a],b)=0$
\end{center}
becomes $\omega ([a,b],i)=0$ and thus $[a,b]\in I^{{\bot}}$. Furthermore as
\begin{center}
$\omega([a,b],c)=0$ for $c\in I^{{\bot}}$
\end{center}
we have $[a,b]\in I$. Hence $K$ is normal in $H^{{\bot}}$.

On the other hand, it is clear that $R$ is connected and
$dim(R)<dim(H^{{\bot}})$. Also, $R$ is a reduced symplectic manifold
of $(G,\omega^{+})$ because $I\cap I^{{\bot}}$ is the nilradical
of the restriction of $\omega$ to $I^{{\bot}}$.

From the identity for $b\in I\cap I^{{\bot}},x\in I^{{\bot}}$ and
$ z\in I$ (or $ z\in I^{{\bot}}$),
\begin{center}
$\omega(bx,z)=-\omega(x,[b,z])=0$
\end{center}
it follows that $I\cap I^{{\bot}}$ is a 2-sided ideal of $I^{{\bot}}$.
All this implies that (3) is an exact sequence of affine Lie groups.

Finally the formula
\begin{center}
$\omega(bx,a)=-\omega(x,[b,a])=0$
\end{center}
where $ b,x\in I\cap I^{{\bot}}$ and $a\in L(G)$ shows that $bx=0$. This means that the affine structure of $K$ deduced from $\nabla^{+}$ is trivial. \qed

\begin{theorem}\label{theo2.4} Let $(G,\omega^{+})$ be a connected symplectic Lie group and $H$ a non-trivial connected normal Lie subgroup of $G$. If $L_{H}$ is Hamiltonian
on $G$ then $H^{{\bot}}$ is an affine closed subgroup of $(G,\nabla^{+})$,
the homogeneous manifold $G/H^{{\bot}}$ is endowed with a commutative parallelism
(and consequently with an affine structure) and the canonical maps in
\begin{equation}
H^{{\bot}}\rightarrow G {\rightarrow} G/H^{{\bot}}
\label{eq:ED4}
\end{equation}
are affine maps.\end{theorem}

\proof Since $L_{H}$ is Hamiltonian, the leaves of ${\cal F^{{\bot }}}$
are closed. In particular $H^{{\bot}}$ is closed in $G$ and (4)
is a principal bundle whose fiber is the affine Lie group $H^{{\bot}}$.
Let $\{x_{1},..., x_{r}\}$ be a basis of $I=L(H)$. We know that the foliation
${\cal F^{{\bot }}}$  is given by the intersection of the kernels of
the left-invariant forms $\eta_{j}^{+}:=i(x_{j}^{+})\omega^{+}$
but also by the intersection of the closed 1- forms  $\eta_{j}:=i(x_{j}^{-})\omega^{+}$
with $1\leq j\leq r$. The 1-forms $\eta_{j}$ are exact and basic relatively to
the projection $p$. So,the projected 1-forms $\eta_{j}$ define a
local parallelism on $G/H^{{\bot}}$. But this parallelism is global and
commutative because the $\eta_{j}$ are exact. Hence the injection $i$
of $H^{{\bot}}$ in $G$ and the projection $p$ in (4) are affine maps.\qed

\section{Momentum maps in Symplectic Lie groups}
Let $(G,\omega^{+})$ be a connected symplectic Lie group. The scalar 2-cocycle
$\omega $ determines a representation of the Lie algebra ${\cal G}$ by endomorphisms of the affine space ${\cal G}^{*}$ given by $x\rightarrow (i(x)\omega ,ad^{*}(x))$, where
$x\rightarrow ad^{*}(x))$ is the coadjoint representation of ${\cal G}$.
This affine representation provides a homomorphism of Lie groups of $\tilde G$,
the universal covering group of $G$, into the classical affine group
$Aff({\cal G}^{*})\simeq  {\cal G^{*}}\times GL({\cal G^{*}})$ (semi-direct product) given
for $\sigma =exp(x)$ with $x\in {\cal G}$ by the formula
\begin{equation}
\rho (\sigma )=(\tilde Q(\sigma ),Ad^{*}(\sigma ))
\label{eq:ED5}
\end{equation}
where
\begin{equation}
\tilde Q(\sigma) =\Sigma _{k=1}^\infty (1/k!)(ad^{*}(x))^{k-1}.i(x)\omega
\label{eq:ED6}
\end{equation}

Notice that $\tilde G$ is not necessarily a exponential Lie group.

As $\omega$ is non-degenerate, the orbit $O$ of $0\in {\cal G}$ by the $\rho $ action,
is open. Let $\nabla^{0}$ be  the connection induced in $O$ by the usual connection
of the affine space ${\cal G}^*$. The pullback of $\nabla^{0}$ by the orbital map
from  $\tilde G$ in $O$ is a left-invariant affine structure. If $\Pi$
is the covering map from $\tilde G$ to $G$, there exists a unique
left-invariant affine structure on $G$ such that $\Pi$ is affine. This latter affine structure on $G$ coincides with that given by formula (\ref{eq:ED1}). Moreover, there exists a representation $h$ of the fundamental group $\Pi_1(G)$ of $G$ by affine transformations of ${\cal G}^{*}$ such that $h(\gamma )\circ \tilde Q=\tilde Q\circ \gamma $ for
every $\gamma $ in $\Pi_1(G)$.

The map $\tilde Q$ is a developing map and $h$ is the holonomy representation
of the affine structure $\nabla^{+}$. The geodesical completeness of $\nabla^{+}$ is equivalent to the fact that $\tilde Q$ is a diffeomorphism.

Suppose that $J:G\rightarrow {\cal G}^{*}$ is a momentum map
for $L_{G}$. This is case if $H^{1}(G,\mathbb{R})=0$ or $H^{1}({\cal G}^{*})=H^{2}({\cal G}^{*})=0$, for instance. This means that, for every $x\in {\cal G}$, the derivative
of the map $J_{x}:G\rightarrow {I\!\!R} $ with $J_{x}(\sigma ):=J(\sigma)(x)$,
is given by $i(x^{-}){\omega}^{+}$. The map $\psi _{(\sigma ,x)}:G\rightarrow {I\!\!R}$
given by
\begin{equation}
\psi _{(\sigma ,x)}(\tau ):=J_{x}(L_{\sigma }(\tau ))-J_{Ad({\sigma}^{-1})(x)}(\tau )
\label{eq:7}
\end{equation}
with $\sigma \in G$ and $x\in {\cal G}$ is constant, because $G$ is
connected. Taking $\tau=\epsilon $ we obtain a mapping  $Q:G\rightarrow {\cal G^{*}}$, that is the 1-cocycle of $G$ relative to
the coadjoint representation, associated with $J$. For $\sigma \in {\cal G}$ putting
\begin {center}
$J'(\sigma ):=J(\sigma )-J(\epsilon )$
\end {center}
defines another momentum map with associated 1-cocycle $Q'=J'$.

We put:

\begin{definition}\label{def3.1} Let $(G,\omega {+})$ be a Hamiltonian
symplectic Lie group. If $J$ is a momentum map for $L_{G}$ with 
$J$ as 1-cocycle we will say that $J$ is a
momentum-cocycle.
\end{definition}

Hence we have

\begin{lemma}\label{lem3.2} Let $(G,\omega {+})$ be a symplectic Lie group. If
$L_{G}$ is Hamiltonian we have:
\\
1. There exists a momentum-cocycle $J$ for $L_{G}$ such that $J(\epsilon )=0\in {\cal G^{*}}$.
\\
2. If $J$ is a momentum-cocycle then the image of $J$ is diffeomorphic
to the manifold $G/D$ where $D:=\{\sigma \in G;J(\sigma)=J(\epsilon )\}$
is a discrete subgroup of $G$. Consequently every momentum map
for $L_{G}$ is a local diffeomorphism.
\end{lemma}

\proof Let $J=Q$ be a momentum-cocycle with $J(\epsilon )=0\in {\cal G}^{*}$.
The map given by $\rho (\sigma ):=(Q(\sigma ),Ad^{*}(\sigma ))$ with
$\sigma \in G$, is a representation of $G$ by affine transformations of
${\cal G^{*}}$ such that the orbit $O$ of $0\in {\cal G^{*}}$ is open. 
Hence, $Im J=O$ can be identified with $G/D$ where $D$ is the $\rho $-stabilizer
in $0\in {\cal G^{*}}$. \qed

\begin{remark}\label{rem1} If $L_{G}$ is Hamiltonian there exists a momentum-cocycle
for $L_{G}$ which is the developing map of the affine structure
given by formula (\ref{eq:ED1}).\end{remark}

A direct calculation shows the following
fact
\begin{lemma} \label{lem3.3} If $(\tilde G,\omega ^{+})$ is a simply connected
symplectic Lie group, the map given by formula (\ref{eq:ED6}) is a
momentum-cocycle for $L_{\tilde G}$.\end{lemma}

In what follows we want to find necessary and sufficient conditions 
to have Hamiltonian symplectic Lie groups.

Let $J:G\rightarrow {\cal G}$ be a momentum map for $L_{G}$ and
$\Pi:\tilde G\rightarrow G$ the universal covering map. A simple
calculation shows that $\tilde {J}:=J\circ\Pi$ is a momentum map
for $L_{\tilde G}$. Moreover if we suppose that $J(\epsilon)=0$,
it follows that
\begin{equation}
\Pi^{-1}(\epsilon)=:\Pi_{1}(G)\subset {\tilde J}^{-1}(\{0\})
\label{eq:8}
\end{equation}
We have also

\begin{lemma} \label{lem3.4} Let $\tilde{J}$ be a momentum-cocycle for
$L_{\tilde{G}}$ with $\tilde{J}(\tilde{\epsilon})=0\in {\cal G^{*}}$.
Then $\tilde{J}$ factors, by means of $\Pi$, into a momentum map
$J$ for $L_{G}$ if and only if condition (8) is verified.\end{lemma}

\proof Put $J(\sigma):= \tilde{J}(\tilde{\sigma})$ where
${\tilde {\sigma}}\in\Pi^{-1}(\sigma)$. We must show that this
expression defines a map.In other words, we must verify that
we have $ \tilde{J}(\tilde{\sigma})=\tilde{J}(\tilde{\tau})$
for every  $\tilde{\tau}\in\Pi^{-1}(\sigma)$. It is clear that
it suffices to consider the case where $\sigma=\epsilon$.
As $\tilde {\sigma}, \tilde {\tau}\in\Pi^{-1}(\epsilon)=\Pi_{1}(G)$
it follows that $\tilde {\sigma}^{-1}\tilde {\tau}\in\Pi^{-1}(\epsilon)$
and hence
$\tilde {\sigma}^{-1}\tilde {\tau}\in  \tilde{J}^{-1}(\{0\})$.
Since $\tilde{J}$ is a 1-cocycle,
$\tilde{J}(\tilde {\sigma}^{-1}\tilde {\tau})=0$
can be written as
\begin {center}
$\tilde{J}(\tilde {\sigma}^{-1})+Ad_{\tilde{G}}^{*}(\tilde {\sigma}^{-1})(\tilde{J}(\tilde {\tau}))=0$.
\end {center}
Applying $Ad_{\tilde{G}}^{*}(\tilde {\sigma})$ to this equality
we get
\begin {center}
$Ad_{\tilde{G}}^{*}(\tilde {\sigma})(\tilde{J}(\tilde {\sigma}^{-1}))+\tilde{J}(\tilde {\tau})=0$.
\end {center}
But the first term on the left hand side is equal to  $-\tilde{J}(\tilde {\sigma})$, consequently we have
 $ \tilde{J}(\tilde{\sigma})=\tilde{J}(\tilde{\tau})$
and the map $J$ is well defined.
Finally a simple calculation shows that $J$ is a
momentum-map for $L_{G}$. \qed

These preliminary facts imply the following important result

\begin{theorem}\label{theo3.5} Let $(G,\omega^{+})$ be a  connected
symplectic Lie group. Then, $L_{G}$ is Hamiltonian
if and only if the holonomy representation of the
affine Lie group $(G,\nabla^{+})$, is trivial.\end{theorem}

\proof Suppose $L_{G}$ Hamiltonian. There exists
a momentum map $J$ for $L_{G}$ such that $J(\epsilon)=0$.
If we put $\tilde{J}:=J\circ\Pi$ we get a
momentum-cocycle for $L_{\tilde G}$ such that
$\tilde {J}(\epsilon)=0$ and (8) is verified.
So, $h(\gamma)=id_{\cal G^{*}}$ for every $\gamma\in\Pi_{1}(G)$.
This means that for each $\tilde{\sigma}\in\tilde G$,
\begin {center}
$\tilde{J}(\gamma\tilde{\sigma})=\tilde{J}(\tilde{\sigma})$
\end {center}
In fact we have
\begin {center}
$\tilde{J}(\gamma\tilde{\sigma})=\tilde{J}(\gamma)+Ad_{\tilde{G}}^{*}( \gamma)(\tilde{J}(\tilde \sigma))=\tilde{J}(\gamma)+\tilde{J}(\tilde{\sigma})
=\tilde{J}(\tilde{\sigma})$
\end {center}

Reciprocally, if $h(\gamma )=id_{\cal G^{*}}$ for every $\gamma \in\Pi _{1}(G)$,
it follows that we have for $\tilde\tau=\gamma \tilde\sigma $,
\begin {center}
$\tilde{J}(\tilde{\tau})=\tilde{J}(\gamma\tilde{\sigma})=h(\gamma )\tilde{J}(\tilde{\sigma})=\tilde{J}(\tilde{\sigma})$.

\end {center}

Hence, we get a momentum-cocycle for $L_{G}$ by setting $J(\sigma ):=\tilde{J}(\tilde{\sigma})$
where $\tilde{\sigma}\in\Pi ^{-1}(\sigma )$ and $\tilde{J}$ is a momentum-cocycle
for $L_{\tilde {G}}$. \qed

\begin{remark}\label{rem1} Let $(T^{2n},\omega ^{+})$ be a symplectic torus. That this group is not
Hamiltonian follows from the compactness of $ T^{2n}$ or from the fact that the
holonomy of $ (T^{2n},\nabla^{+})$ is not trivial.\end{remark}

For the sake of completeness let us recall at this point that Theorem 3 in \cite{Li-Me}  can
be presented as follows

\begin{theorem}\label{theo3.6} Let $(G,\omega^{+})$ be a connected symplectic Lie group,
$\Pi: \tilde {G}\rightarrow G$ the universal covering of $G$ and
$J$ a momentum-cocycle for $L_{\tilde G}$. Then, $J$ is a global diffeomorphism
i.e. $\nabla^{+}$ is geodesically complete, if and only if $G$ is unimodular.
Moreover, if this condition is verified  then $G$ is solvable.

\end{theorem}

\begin{example}Consider the manifold $G=\mathbb{R_{+}^{*}}\times\mathbb{R} $ endowed with the product $(a,b)(c,d):=(ac,ad+b)$.The corresponding Lie algebra ${\cal G}=span\{e_{1},e_{2}\}$ with $[e_{1},e_{2}]=e_{2}$ endowed with the scalar 2-cocycle $\omega=e^{*}_{1}\wedge^{*}_{2}$ is a  symplectic Lie algebra, where $\{e^{*}_{1},e^{*}_{2}\}$ denotes the dual basis. As $G$ is simply connected it is a Hamiltonian Lie group. We will determine a moment-cocycle for $L_{G}$ using the formula (6).

For $x=x_{1}e_{1}+x_{2}e_{2}\in {\cal G}$ and $t\in \mathbb{R}$ a direct calculation gives:
\begin{center}
$exp(tx) = (1,tx_{2})$ if $x_{1}=0$ ;
$exp(tx)= (e^{x_{1}},\frac{x_{2}}{x_{1}}(e^{x_{1}t}-1))$ if $x_{1}\neq 0$ 

$ad^{*}_{x}(e_{1})=0 $ ; $ad^{*}_{x}(e_{2})=x_{2}e^{*}_{1}-x_{1}e^{*}_{2}$
\end{center}

This implies 
\begin{center}
$Q(exp x)= -x_{2}e_{1}^{*} , x_{1}=0$ 

$Q(exp x)= \frac{x_{2}}{x_{1}}(e^{-x_{1}}-1)e^{*}_{1}+(1-e^{-x_{{1}}})e^{*}_{2} , x_{1}\neq 0$
\end{center}

Hence for $(a, b)\in G $ we get then
\begin{center}
$Q(1,b) = -be^{*}_{1}$

$Q(a,b)=\frac{b}{a-1}(\frac{1}{a})e^{*}_{1}+(1-\frac{1}{a})e^{*}_{2}$ ,  otherwise
\end{center}
Recall that $Q$ is also the developing map of the canonical affine structure. Taking $H= \{(1,b); b\in \mathbb{R}\}$, as in Theorem 1, we get $K=R=\{0\}$.
Moreover $H$ is an affine Lie subgroup of $G$ and the canonical sequence, $H\rightarrow G\rightarrow \mathbb{R}e_{1}$ is an exact sequence of affine Lie groups.
Also $({\cal G},\omega,j)$ , with $j(e_{1})= e_{2}$ and $j(e_{2})=- e_{1}$, is a positive Kahlerian Lie algebra non Ricci flat (see [10]).
 
\end{example}

\begin{example}
Let ${\cal G}_{1}:=Span \{e_{1},e_{2}\}$ the real Abelian Lie algebra endowed with the scalar 2-cocycle $\omega_{1}( e_{1},e_{2})=-1$ and the complex structure $j_{1}(e_{1})=e_{2}$ , $j_{1}(e_{2})=-e_{1}$. As $g_{1}(x,x)=\omega (j(x),x) $ for $x\in {\cal G}_{1}$ is an Euclidian structure, the triple $({\cal G}_{1},\omega_{1},j_{1})$ is a Kahler Lie algebra ( [5]). Denote by $({\cal G}_{2},\omega_{2},j_{2})$ another copy of the same algebra with ${\cal G}_{2}=Span \{e_{3},e_{4}\}$ and by  $\eta$ the representation of  ${\cal G}_{2}$ in $so({\cal G}_{1})$ given by $\eta(e_{3})(e_{1})=e_{2}$, $\eta(e_{3})(e_{2})=-e_{1}$, $\eta(e_{4})(e_{1})=-e_{2}$, $\eta(e_{4})(e_{2})=e_{1}$.
The Lie algebra ${\cal G}$, semi-direct product of ${\cal G}_{1}$ by ${\cal G}_{2}$ via $ \eta$, becomes a positive Kahlerian Lie algebra taking $\omega= \omega_{1}+ \omega_{2}$ and $j=j_{1}+ j_{2}$. The Riemann metric is $g=g_{1}+ g_{2}$ ; this metric is flat. Hence  a priori, on the corresponding simply connected Lie groupe groupe $G=G_{1}\times G_{2}= \mathbb{R}^{4}$ there are two left invariant affine structures and consequently two developing maps and $(G,\omega^{+})$ is Hamiltonian.  It is clear that the orthogonal symplectic of $H:=G_{1} $ is $G_{2}$ and a direct calculation shows that the product in $G$ for $\sigma =(a,b,c,d)$ and $\tau=(a',b',c',d')$ is given by,

\begin{center}

$\sigma\tau= (a+a'cos(c-d)-b'sin(c-d), b+a'sin(c-d)+b'cos(c-d), c+c',d+d')$

\end{center}
As $G$ is not a exponential Lie group we can using the closed $1$-formes on $G$ given by $\omega^{+}(x^{-}, )$ for find the momentum map.  

\end{example}

\section{Structure of simply connected nilpotent symplectic Lie groups}

Let $(G,\omega^{+})$ be a simply connected nilpotent symplectic Lie group.
Consider the exponential Lie subgroup $H:=\{exp tz;t\in I\!\!R\}$
where $z$ is a non-zero central element of ${\cal G}$. Since the
cohomology group $H^{1}(G,I\!\!R)$ is trivial, the action
$L_{G}$ is Hamiltonian. Obviously $J_{H}:=i^{*}\circ J$ is a momentum map for
$L_{H}$ and $J_{H}(exp x)$ is the linear map on $I\!\!Rz$ given by
$z\rightarrow \omega (x,z)$. Hence $J_{H}^{-1}(\{0\})=:H^{\bot }$
is a closed and connected subgroup of $G$ which is the (connected)
 subgroup of $G$  orthogonal to $H$ relatively to $\omega^{+} $. 
Applying Theorem \ref{theo2.3} and Theorem \ref{theo2.4} we obtain the exact canonical
sequence of affine groups:
\begin{equation}
\{\epsilon \}\rightarrow H\rightarrow H^{\bot }\rightarrow R:=H^{\bot }/H\rightarrow \{\epsilon \}
\label{eq:9}
\end{equation}
and the canonical principal fiber bundle
\begin{equation}
\{\epsilon \}\rightarrow H^{\bot}\rightarrow G\rightarrow G/H^{\bot }
\label{eq:10}
\end{equation}
where the projection can be identified with the momentum map $J_{H}$.

But as $z$ is central we can deduce that the Lie algebra of $H^{\bot}$
is an ideal of the Lie algebra ${\cal G}$. As a consequence $G/H^{\bot }$
is a (connected) one-dimensional affine subgroup of $(G,\nabla ^{+})$.

In conclusion we have proved the following result:

\begin{theorem}\label{theo4.1} Let $(G,\omega^{+})$ be a simply connected
nilpotent symplectic Lie group. If $H$ is an exponential 1-dimensional
central Lie subgroup of $G$ then $H$ and $H^{\bot }$ are affine
subgroups of $(G,\nabla^{+})$ and $G$ is a semi-direct product of the normal subgroup $H^{\bot}$
by the group $G/H^{\bot }$ where $H^{\bot}$ is a central extension
of the reduced symplectic Lie group $R=H^{\bot}/H$.\end{theorem}

\begin{remark}\label{rem1} At the level of Lie algebras, Theorem \ref{theo4.1}
corresponds to a particular case of the notion of symplectic double extension
of symplectic Lie algebras introduced by Medina and Revoy in \cite{Me-Re} (see also
\cite{Da-Me1}, \cite{Da-Me2}).\end{remark}

 For Lie groups we put

\begin{definition}\label{def4.2} A symplectic Lie group $(G,\omega ^{+})$ is called
a symplectic double extension of a symplectic Lie group $(R,\eta ^{+})$
if $(G,\omega ^{+})$ contains a totally isotropic Lie subgroup $H$
such that $H^{\bot}/H$ is a symplectic Lie group isomorphic to
$(R,\eta ^{+})$.\end{definition}

In these terms we can say that every simply connected nilpotent symplectic
Lie group $G$ of dimension $2(n+1)$ is a symplectic double extension of
a $2n$-dimensional simply connected nilpotent symplectic Lie group.
Moreover G is obtained by a sequence of symplectic double extensions
starting from the simply connected Abelian Lie group.

\section{The classic affine group as Frobenius Lie group}

\begin{definition}\label{def5.1} A symplectic Lie group $(G,\omega ^{+})$ is called
a Frobenius Lie group if $\omega ^{+}=d\nu ^{+}$ where $\nu^{+}$ is
a left-invariant 1-form on $G$.\end{definition}

Every Frobenius Lie group  $(G,d\nu ^{+})$ is a Hamiltonian Lie group.
In fact the map defined by
\begin{equation}
J(\sigma)(x):=-\nu^{+}_{\sigma }(x_{\sigma }^{-})
\label{eq:11}
\end{equation}
for $\sigma\in G$ and $x\in {\cal G}$ is a $Ad_{G}^{*}$-equivariant
momentum map and $J'(\sigma ):=J(\sigma )-J(\epsilon )$
is a momentum-cocycle for $L_{G}$.

Consequently if $H$ is a non trivial connected normal and
closed subgroup of $(G,d\nu ^{+})$ such that $K:=H^{\bot}\bigcap H$
is non-discrete, then $G$ is a symplectic double extension of
the symplectic Lie group $R:=H^{\bot}/K$. The case where
$K=H$ is particularly nice and simple. To illustrate this situation
we will study the group of the affine transformations of the
real (or complex) affine space of dimension n.

Consider the classical affine group $G=GA(I\!\!R^{n})\equiv I\!\!R^{n}\times GL(I\!\!R^{n})$
(semi-direct product by means of the trivial representation) and let ${\cal G}\equiv I\!\!R^{n}\times gl(I\!\!R^{n})$ be its
Lie algebra. To any $\alpha\in{L(G)^{*}}$ we can associate a unique
pair $(g,M)$ with $g\in ( I\!\!R^{n})^{*}$ and $M\in gl(I\!\!R^{n})$
given by:

\begin {center}
$\alpha (x,u)=g(x)+tr(M\circ u)$
\end {center}
for $x\in I\!\!R^{n}$ and $u\in gl(I\!\!R^{n})$.

The associated 2-coboundary $\delta \alpha $ can be written
\begin {center}
$\delta \alpha ((x,a);(y,v))=-g(u(x))-v(y)-tr(M\circ [u,v])$.
\end {center}

Suppose that the coadjoint orbit of $\alpha\equiv(g,M)$ is open (see Theorem 2.5 of \cite{B-M-O})
and take $\omega^{+}=\delta\alpha^{+}$ and $H=I\!\!R^{n}$.
The momentum map $J:G\rightarrow{\cal G^{*}}$,
$J(\sigma):=Ad_{\sigma}^{*}(\alpha)$ can be written

\begin {center}
$J(\sigma)(x)=\alpha_{\sigma}^{+}(x_{\sigma}^{-})$
\end {center}
where $x^{-}$ denotes the right-invariant vector
field on $G$ associated with $x\in{\cal G}$.

Remember that $J_{H}(\sigma)$ is given by
the restriction of $J(\sigma)$
to the Lie algebra of $H$.

Dardi\'e and Medina have described the architecture of 
the classical affine group in the following terms (result
not yet published)

\begin{theorem}\label{theo5.2} If the classical affine Lie group
$G=GA(I\!\!R^{n})$ is endowed with the left-invariant
symplectic form $\omega^{+}=\delta\alpha^{+}$ with
$\alpha\equiv(g,M)$ and $H=I\!\!R^{n},$ we have:

1. $J_{H}^{-1}(g)$ is a closed subgroup of $G$
containing $H$.

2. The exact canonical sequence of Lie groups
\begin {center}
$\{\epsilon\}\rightarrow H\rightarrow J_{H}^{-1}(g)\rightarrow J_{H}^{-1}(g)/H \rightarrow\{\epsilon\}$
\end {center}
splits and is a sequence of affine Lie groups.

3. The symplectic reduced Lie group $J_{H}^{-1}(g)/H$ is
isomorphic to $GA(I\!\!R^{n-1})$.

4. In the canonical fiber bundle
\begin {center}
$J_{H}^{-1}(g)\rightarrow GA(I\!\!R^{n})\rightarrow (I\!\!R^{n})^{*}-\{0\}$
\end {center}
the fiber is an affine Lie subgroup of $GA(I\!\!R^{n})$
and the projection is affine relatively
to the usual affine structure of $(I\!\!R^{n})^{*}-\{0\}.$
\end{theorem}

The proof of this Theorem is rather technical. Basically it requires
a good understanding of the ideas and tools
developed in the above sections and the following technical Lemma.

\begin{lemma}\label{lem5.3} For $G$ and $H$ as in Theorem \ref{theo5.2} the momentum map
$J_{H}$ is a surjective submersion on $(I\!\!R^{n})^{*}-\{0\}$ and
\begin{equation}
J_{H}(x,T)={}^{t}T^{-1}g
\label{eq:12}
\end{equation}
for $(x,T)\in G$.\end{lemma}

\proof The Lemma is a consequence of the following formulas
describing the coadjoint representation of $G$:
$$Ad^{*}_{(0,U)}(h,N)=({}^{t} U^{-1}(h),U\circ N\circ U^{-1})$$
$$Ad^{*}_{(x,I)}(h,N)=(h,N+h\times x)$$
where $U$ is taken in the connected component in $GL(n,I\!\!R)$ of the unit
$I$ and $(h\times x)(y):= h(y)x$ for $x,y \in I\!\!R^{n}$.

Let us now sketch the proof of the Theorem. Formula (12)
implies that $J_{H}^{-1}(g)=\{(x,T)\in G,{}^{t}T^{-1}g=g\}$ is a
closed subgroup of $G$. The factor group $J_{H}^{-1}(g)/H$ can be
identified with subgroup $\{(0,T)\in G,{}^{t}T^{-1}g=g\}$ of $G$
and consequently it can be interpreted as the group of affine
transformations of $Ker(g)$. On the other hand we have
$L(J_{H}^{-1}(g))=L(H)^{\bot}$. Hence assertions 2. and 4.
are consequences of Theorems \ref{theo2.3} and \ref{theo2.4}. Finally  Lemma \ref{lem5.3} 
implies assertion 3 and Theorem \ref{theo5.2} is proved. \qed

Denote now by $\alpha_{1}$ and $\alpha_{2}$ the elements of
${\cal G}^*$ given by $\alpha_{1}(x,T)=g(x)$ and $\alpha_{2}(x,T)=tr(M{\circ}T)$.
We have $\delta\alpha=\delta \alpha_{1}+\delta\alpha_{2}$ and ${\cal G}$ is the orthogonal direct sum of the symplectic subalgebras ${\cal G}_{1}:=Rad(\delta\alpha_{2})$ and ${\cal G}_{2}:=Rad(\delta \alpha_{1})$(see \cite{B-M-O}). Moreover ${\cal G}_{2}$ is isomorphic to the Lie algebra $aff(Ker(g))$ and ${\cal G}_{1}$ is isomorphic to the natural semidirect product
of $\mathbb R^{n}$ by the commutant of $M$ in $gl(\mathbb R^{n})$ denoted in the following by $C(M)$.

Consequently the infinitesimal version of Theorem \ref{theo5.2}  is given by the following assertion

\begin{proposition}\label{prop5.4} In the notations of Theorem  \ref{theo5.2} if $I=L(H)$,the canonical sequence
\begin {center}
$0\rightarrow I\rightarrow I^{\bot}\rightarrow I^{\bot}/I\rightarrow 0$
\end{center}
is a split exact sequence of Lie algebras with $ I{\cong}I\!\!R^{n}$,$I^{\bot}=I\times {\cal G}_{2}{\cong}\mathbb R^{n}\times aff(\mathbb R^{n-1})$ and $I^{\bot}/I{\cong}{\cal G}_{2}{\cong}Aff(\mathbb R^{n-1})$. It is also an exact sequence of left-symmetric algebras.

Moreover the vector space $aff(I\!\!R^{n})$ is a direct sum of the Lie subalgebras $I^{\bot}$
and $C(M)$.\end{proposition}

Let us be more precise about the  morphisms of Lie algebras in Proposition \ref{prop5.4}. The isomorphism between
${\cal G}_{2}$ and $ aff(I\!\!R^{n-1})$ is given by the choice of a supplementary of $Ker(g)$
into $I\!\!R^{n}$. In fact if $x\in I\!\!R^{n}$ verifies $g\circ M^{i}x=0$ for $0\leq i\leq n-2$
and $g\circ M^{n-1}x=1$, then the map $(0,u)\rightarrow(u(M^{n-1}x),u')$, where $u'$ is the restriction of $u$ to $Ker(g)$ ,defines such an isomorphism. Moreover the symplectic reduced form
corresponds to the 2-coboundary associated with $(g_{1},M_{1})\in aff((I\!\!R^{n-1})^{*})$ where
$g_{1}$ and $M_{1}$ are given by
$$tr(M\circ u)=g_{1}(u(M^{n-1}(x)))+tr(M_{1}\circ u')$$

Suppose $M$ is nilpotent (this is always possible, see\cite{B-M-O}), then $g_{1}$ can be identified with
${}^tM(g)\in (\mathbb R^{n})^{*}$ and $M_{1}$ is also nilpotent. 

Hence in the case where $\alpha\equiv(g,N)$ with $N$ nilpotent, an iteration of the (infinitesimal) procedure of symplectic reduction yields a decomposition of the Lie algebra $aff(\mathbb R^{n})$ as a sum of Abelian subalgebras,

$$aff(\mathbb R^{n})=K_{n}\oplus K_{n-1}\oplus ... \oplus K_{1}\oplus C(N)\oplus C(N_{1})\oplus ...\oplus C(N_{n-1})$$

where $K_{i}$ and $C(N_{i-1})$ are the subalgebras corresponding respectively to the (totally) isotropic ideal and the commutant of $N_{i-1}$ found at the moment of the i-reduction.

Let us explain in few words this decomposition.

The elements of $C(N_{i})$ are the polynomials in $N_{i}$ with $N_{i}$ the endomorphism defined by $N^{i}_{i}=0$ and $N^{j}_{i}x=N^{j}x$ for $j<i$ and $0\leq i\leq {n-1}$.

Let $B=\{x,Nx,...,N^{n-1}x\}$ be the dual basis of 
$$
B^{*}=\{{}tN^{n-1}g,{}^tN^{n-2}g,...,g\}.
$$
Using $B$ we can identify the Lie algebra
$aff(\mathbb R^{n})$ with the Lie algebra of $(n+1)\times(n+1)$-matrices with the last row equal to zero.

In this context the coefficients of the $n\times n$ linear part $(a_{{i},{j}})$ of any element of $K_{i}$ are null 
except for $a_{{1},{i}},a_{{2},{i}},...a_{{i-1},{i}}$ of the $i-column$.


Consequently, relatively to the basis $B$ dual of $B^{*}$, the subspace $L=K_{n}\oplus K_{n-1}\oplus...\oplus K_{1}$ can be identified with the Lie subalgebra of $gl(\mathbb R^{n+1})$ consisting of strict upper triangular matrices.

Finally, in the basis $B$ given above, $L'=C(N)\oplus C(N_{1})\oplus...\oplus C(N_{n-1})$ could be identified with the subalgebra of $(n+1)\times (n+1)$ lower triangular matrices having the last row equal to zero. 

As at each symplectic reduction $K_{i}$ and $C(N_{i})$ are totally isotropic subspaces, it follows that $L$ and $L'$ are Lagrangian subalgebras of $aff(\mathbb R^{n})$ relatively to $\delta\alpha$ with $\alpha\equiv(g,N)$ and $N$ nilpotent.

Denote by $\Lambda$ and $\Lambda'$ the connected Lie subgroups of $Aff(\mathbb R^{n})$ corresponding to subalgebras $L$ and $L'$. The canonical left actions $L_{\Lambda}$ and
$L_{\Lambda'}$ on $Aff(I\!\!R^{n})$ are Hamiltonian and consequently $\Lambda$ and $\Lambda'$
are closed (Theorem 3.1 of \cite{B-M-O}).

However we have (see \cite{D}):

\begin{proposition}\label{prop5.5}
 i.$Aff(I\!\!R^{n})$ has exactly two open orbits for the coadjoint representation. Moreover,two elements $(g,M)$ and $(h,N)$ having an open orbit are in the same orbit if and only if the bases $\{{}^tM^{n-1}g,{}^tM^{n-2}g,...,g\}$ and $\{{}^tN^{n-1}g,{}^tN^{n-2}g,...,g\}$ give the same orientation.

ii.If $\Omega$ and $\Omega'$ are two left invariant symplectic forms on $Aff(I\!\!R^{n})$, there exists an automorphism $\varphi$ of the Lie algebra $aff(I\!\!R^{n})$ such that
$$\Omega_{\epsilon}(.,.)=\Omega'_{\epsilon}(\varphi.,\varphi.)$$\end{proposition}

Proposition \ref{prop5.4} and the above analysis imply:

\begin{theorem}The symplectic Lie group $(Aff(\mathbb R^{n}),d\alpha^{+})$ is endowed with a
pair of left invariant transversal Lagrangian foliations whose leaves are closed and given by $\Lambda$ and $\Lambda'$.\end{theorem}


\end{document}